\documentclass[12pt,a4paper]{article}
\usepackage{a4,amsmath,amssymb,amsfonts,amsthm}
\usepackage{amsrefs}
\usepackage[small,nohug,heads=vee]{diagrams}
\usepackage[all]{xy}
\usepackage{palatino}

\begin{document}
\bibliographystyle{plain}


\def\mR{\M{R}}
\def\mZ{\M{Z}}
\def\mN{\M{N}}           
\def\mQ{\M{Q}}
\def\mC{\M{C}}
\def\mG{\M{G}}
\def\mF{\M{F}}



\def\Spec{{\rm Spec}}
\def\rg{{\rm rg}}
\def\Hom{{\rm Hom}}
\def\Aut{{\rm Aut}}
 \def\Tr{{\rm Tr}}
 \def\Exp{{\rm Exp}}
 \def\Gal{{\rm Gal}}
 \def\End{{\rm End}}
 \def\det{{{\rm det}}}
 \def\Td{{\rm Td}}
 \def\ch{{\rm ch}}
 \def\che{{\rm ch}_{\rm eq}}
  \def\Spec{{\rm Spec}}
\def\Id{{\rm Id}}
\def\Zar{{\rm Zar}}
\def\Supp{{\rm Supp}}
\def\eq{{\rm eq}}
\def\Ann{{\rm Ann}}
\def\LT{{\rm LT}}
\def\Pic{{\rm Pic}}
\def\rg{{\rm rg}}
\def\et{{\rm et}}
\def\sep{{\rm sep}}
\def\ppcm{{\rm ppcm}}
\def\ord{{\rm ord}}
\def\Gr{{\rm Gr}}
\def\ker{{\rm ker}}
\def\rk{{\rm rk}}
\def\Coh{{\rm Coh}}


\def\beginProof{\par{\bf Proof. }}
 \def\endProof{${\qed}$\par\smallskip}
 \def\pr{^{\prime}}
 \def\prpr{^{\prime\prime}}
 \def\mtr#1{\overline{#1}}
 \def\ra{\rightarrow}
 \def\mfp{{\mathfrak p}}

 \def\mQ{{\Bbb Q}}
 \def\mR{{\Bbb R}}
 \def\mZ{{\Bbb Z}}
 \def\mC{{\Bbb C}}
 \def\mN{{\Bbb N}}
 \def\mF{{\Bbb F}}
 \def\mA{{\Bbb A}}
  \def\mG{{\Bbb G}}
  \def\mP{{\Bbb P}}
 \def\CI{{\cal I}}
 \def\CE{{\cal E}}
  \def\CF{{\cal F}}
 \def\CJ{{\cal J}}
 \def\CH{{\cal H}}
 \def\CO{{\cal O}}
 \def\CA{{\cal A}}
 \def\CB{{\cal B}}
 \def\CC{{\cal C}}
 \def\CD{{\cal D}}
 \def\CK{{\cal K}}
 \def\CL{{\cal L}}
 \def\CI{{\cal I}}
 \def\CM{{\cal M}}
  \def\CN{{\cal N}}
\def\CP{{\cal P}}
\def\CR{{\cal R}}
\def\CQ{{\cal Q}}
\def\CG{{\cal G}}
\def\CT{{\cal G}}
 \def\wt#1{\widetilde{#1}}
 \def\mod{{\rm mod\ }}
 \def\refeq#1{(\ref{#1})}
 \def\blb{{\big(}}
 \def\brb{{\big)}}
\def\mc{{{\mathfrak c}}}
\def\mcpr{{{\mathfrak c}'}}
\def\mcprpr{{{\mathfrak c}''}}
\def\ss{{\rm ss}}
\def\parf{{\rm parf}}
\def\P1{{{\bf P}^1}}
\def\cod{{\rm cod}}
\def\pr{\prime}
\def\prpr{\prime\prime}
\def\ss{\scriptstyle}
\def\OX{{ {\cal O}_X}}
\def\mpartial{{\mtr{\partial}}}
\def\inv{{\rm inv}}
\def\indlim{\underrightarrow{\lim}}
\def\prolim{\underleftarrow{\lim}}
\def\pprolim{'\prolim'}
\def\Pro{{\rm Pro}}
\def\Ind{{\rm Ind}}
\def\Ens{{\rm Ens}}
\def\without{\backslash}
\def\pbdb{{\Pro_b\ D^-_c}}
\def\qc{{\rm qc}}
\def\Com{{\rm Com}}
\def\an{{\rm an}}
\def\gfield{{\rm\bf k}}
\def\s{{\rm s}}
\def\dR{{\rm dR}}
\def\ari#1{\widehat{#1}}
\def\ul#1{\underline{#1}}
\def\sul#1{\underline{\scriptsize #1}}
\def\mou{{\mathfrak u}}
\def\ich{\mathfrak{ch}}
\def\cl{{\rm cl}}
\def\K{{\rm K}}
\def\R{{\rm R}}
\def\F{{\rm F}}
\def\L{{\rm L}}
\def\pgcd{{\rm pgcd}}
\def\rc{{\rm c}}
\def\N{{\rm N}}
\def\E{{\rm E}}
\def\H{{\rm H}}
\def\CHOW{{\rm CH}}
\def\A{{\rm A}}
\def\d{{\rm d}}
\def\Res{{\rm  Res}}
\def\GL{{\rm GL}}
\def\Alb{{\rm Alb}}
\def\alb{{\rm alb}}
\def\Hdg{{\rm Hdg}}
\def\Num{{\rm Num}}
\def\Irr{{\rm Irr}}
\def\Frac{{\rm Frac}}
\def\Sym{{\rm Sym}}
\def\indlim{\underrightarrow{\lim}}
\def\prolim{\underleftarrow{\lim}}
\def\red{{\rm red}}
\def\naive{{\rm naive}}
\def\ch{{\rm ch}}
\def\Td{{\rm Td}}
\def\T{{\rm T}}
\def\min{{\rm min}}
\def\slope{{\rm slope}}
\def\max{{\rm max}}
\def\min{{\rm min}}
\def\Sup{{\rm Sup}}
\def\Qb{\bar{\mQ}}
\def\mn{{\mu_n}}
\def\m2{{\mu_2}}
\def\TW{\{-1\}}
\def\BC{{\widetilde{\ch}}}
\def\ol#1{\overline{#1}}
\def\LocF{{\rm LocF}}
\def\LL{{\rm LL}}
\def\Coh{{\rm Coh}}
\def\Perf{{\rm Perf}}


\def\RHom{{\rm RHom}}
\def\rRHom{{\mathcal RHom}}
\def\rHom{{\mathcal Hom}}
\def\dotimes{{\overline{\otimes}}}
\def\Ext{{\rm Ext}}
\def\rExt{{\mathcal Ext}}
\def\Tor{{\rm Tor}}
\def\rTor{{\mathcal Tor}}
\def\SP{{\mathfrak S}}
\def\perf{{\rm perf}}
\def\Bl{{\rm Bl}}

\def\H{{\rm H}}
\def\D{{\rm D}}
\def\Del{{\mathfrak D}}
\def\Stab{{\rm Stab}}
\def\Div{{\rm Div}}
\def\Ver{{\rm Ver}}
\def\insep{{\rm insep}}

 \newtheorem{theor}{Theorem}[section]
 \newtheorem{prop}[theor]{Proposition}
 \newtheorem{propdef}[theor]{Proposition-Definition}
 \newtheorem{cor}[theor]{Corollary}
 \newtheorem{lemma}[theor]{Lemma}
 \newtheorem{sublem}[theor]{sub-lemma}
 \newtheorem{defin}[theor]{Definition}
 \newtheorem{conj}[theor]{Conjecture}
\newtheorem{rem}[theor]{Remark}

 \parindent=0pt
 \parskip=5pt

\author{Damian R\"ossler}

\title{Canonical isomorphisms of determinant line bundles}

\maketitle

\begin{abstract}
We prove a local refinement of the Grothendieck-Riemann-Roch theorem in degree one.
\end{abstract}

\section{Introduction}

The aim of this text is to provide a proof of the following theorem.

Let   $S$ be a  noetherian scheme and let 
$g:Y\to S$ be a smooth and strongly projective morphism. 
Recall that this means that there is an $N\geq 0$ and a factorisation of 
$g$ into a closed $S$-immersion $Y\to\mP^N_S$ followed by projection to $S$. 
Suppose also that $2$ is invertible in $S$. We suppose that $Y$ has constant relative dimension $d$ over $S$.

For any coherent locally free sheaf $F$ on $Y$, we shall write 
$\lambda(F):=\det(\R^\bullet g_*(F))$. Here $\det(\cdot)$ is the Knudsen-Mumford determinant of a perfect complex (note that $\R^\bullet g_*(F)$ is a perfect complex by the semicontinuity theorem because $g$ is proper and flat). We shall denote by 
$\Sym^k(F)$ the $k$-th symmetric power of $F$ and we shall 
write $F^\vee:=\underline{\rm Hom}(F,\CO_X)$ for the dual of $F$. 
If $\CM$ is a line bundle on $Y$ and $k\in\mZ$, we define
$\CM^{\otimes k}:=\otimes_{i=1}^k \CM$ if $k\geq 0$ and 
$\CM^{\otimes k}:=\otimes_{i=1}^{-k} \CM^\vee$ if $k<0$. As is costumary, we shall write 
$\Omega_g=\Omega_{Y/S}$ for the sheaf of differentials of $g$.

Let $\CL$ be a line bundle (ie a locally free sheaf of rank one) on $Y$. 

\begin{theor}
There is a canonical isomorphism
\begin{equation}
\lambda(\CL)^{\otimes 2^{2d+2}}\simeq 
 \bigotimes_{i=0}^{2d}\bigotimes_{j=0}^i 
 \lambda\big(\CL^{\otimes 2}\otimes\Sym^j(\Omega_{Y/S})\big)^{\otimes2^{2d-i}(-1)^j{i\choose j}}.
 \label{mainis}
 \end{equation}
 This isomorphism is compatible with any base change to a noetherian scheme if hypothesis (H) below holds.
 \label{mainth}
\end{theor}
Hypothesis (H) is the assumption (described in section \ref{secducrot} below) that the multiadditivity of Ducrot's generalisation 
of the Deligne pairing is compatible with any base-change to a noetherian scheme. We were not able to verify this in detail (but we believe that it is true) so we prefer to list it as a supplementary assumption.

For example, suppose that $d=1$. We then get an isomorphism 
\begin{equation}
\lambda(\CL)^{16}\simeq\lambda(\CL^{\otimes 2})^{\otimes 7}\otimes\lambda(\CL^{\otimes 2}\otimes\Omega_{Y/S})^{\otimes(-4)}\otimes\lambda(\CL^{\otimes 2}\otimes\Omega_{Y/S}^{\otimes 2})
\label{vardel}
\end{equation}
In particular,  writing $\lambda_k:=\lambda(\Omega_{Y/S}^{\otimes k})$, \refeq{vardel} gives
$$
\lambda_k^{\otimes 16}\simeq\lambda_{2k}^{\otimes 7}\otimes\lambda_{2k+1}^{\otimes(-4)}
\otimes\lambda_{2k+2}.
$$
By Serre duality, there is a canonical isomorphism $\lambda_0\simeq\lambda_1$. Thus, setting $k=0$ we obtain a canonical isomorphism
$$
\lambda_1^{\otimes 13}\simeq\lambda_2.
$$
This is a special case of the Mumford isomorphism (see \cite{Mumford-Stability}). If $Y$ is an elliptic scheme, 
there is a canonical isomorphism $\Omega_{Y/S}\simeq g^*(g_*(\Omega_{Y/S}))$ so that 
we have canonically
$$
\lambda_1\simeq\lambda_k\simeq g_*(\Omega_{Y/S}).
$$
In particular there is a canonical isomorphism $(g_*(\Omega_{Y/S}))^{\otimes 12}\simeq\CO_S$. Possibly up to sign, this 
is the classical discriminant modular form. This suggests that the isomorphism in Theorem \ref{mainth} 
is in some sense optimal.

When $Y$ is an elliptic scheme over $S$ and $\CL$ is a non trivial torsion line bundle, whose order is prime to the characteristic of all the residue fields of $S$, then $\R^\bullet g_*(\CL)=0$. In that case, both sides  
of \refeq{mainis} are canonically isomorphic to the trivial line bundle. Thus the isomorphism \refeq{mainis} provides an element of $\Gamma(S,\CO_S^\ast)$, in other words an elliptic unit. It seems likely that one can construct all the Siegel units in this way but to prove this, one will have probably have to wait for a metric version of Theorem \ref{mainth}. See below for a discussion.

When $d=2$ and $\CL=\CO_Y$, we get the isomorphism 
$$
\lambda(\CO_Y)^{\otimes 64}\simeq 
\lambda(\CO_Y)^{\otimes 31}\otimes\lambda(\Omega_f)^{\otimes(-26)}
\otimes\lambda(\Sym^2(\Omega_f))^{\otimes 16}\otimes
\lambda(\Sym^3(\Omega_f))^{\otimes(-6)}\otimes
\lambda(\Sym^4(\Omega_f)).
$$
from Theorem \ref{mainth}. 
This is equivalent to 
$$
\lambda(\CO_Y)^{\otimes 33}\otimes\lambda(\Omega_f)^{\otimes 26}\otimes
\lambda(\Sym^3(\Omega_f))^{\otimes 6}\simeq\lambda(\Sym^2(\Omega_f))^{\otimes 16}\otimes \lambda(\Sym^4(\Omega_f)).
$$
and there are similar identities in any relative dimension.

Here is our method of proof. We first give a proof of the geometric fixed formula 
for an involution, which avoids any reference to K-theory and uses only 
the geometric properties of quotients. The idea to use quotients to prove the fixed point formula 
is due to  Thomason (see \cite{Thomason-Formule}) and most probably many earlier authors but our proof relies on the crucial fact that when the 
fixed point scheme is a Cartier divisor then the quotient morphism is flat. This seems to 
be well known fact (J. Oesterl\'e kindly explained the proof to me many years ago) but we could find no proof of it in the literature in the required generality and we provide one in Proposition \ref{vimplem} (1). Our proof of the geometric fixed point formula is sufficiently explicit 
to provide canonical isomorphisms at every step (rather than equalities 
in the Picard group) but ends with an error term, which turns out 
to be a line bundle arising from a higher dimensional version of the Deligne 
pairing. This pairing was studied by Ducrot in \cite{Ducrot-Cube} and we use his results 
to show that this line bundle is canonically trivial, compatibly with any base change to a noetherian scheme, conditional on hypothesis (H). 
We then apply this formula to the space $Y\times_S Y$ with the involution swapping the factors. 
Nori (see \cite{Nori-Hirzebruch}) was apparently the first one to notice that the fixed point formula applied to this situation recovers 
the Adams-Riemann-Roch for the Adams operation $\psi^2$ and using our method we thus recover 
a refinement of this formula (in degree one), where isomorphisms of line bundles become canonical and 
the torsion part is controlled uniformly. This is formula \refeq{mainis}.

In \cite{Eriksson-Isomorphisme} Eriksson gives a proof of a functorial refinement of the Adams-Riemann-Roch 
formula (see also \cite{Eriksson-Isomorphisme-CRAS} for an announcement), which can also be used to prove a weaker version of 
Theorem \ref{mainth}. 
It is weaker in the sense that the provided isomorphism, although canonical, will include a $2^\infty$-torsion line bundle, which is undetermined and also because the resulting linear combination in the symmetric powers of $\Omega_{Y/S}$ will a priori depend on the 
dimension of the total space. 

Similarly, using Franke's work in \cite{Franke-RRFF}, it is possible to prove a weak version of Theorem \ref{mainth}, 
where an undetermined (not necessarily $2^\infty$) torsion line bundle will be included 
(but on the other hand the linear combination in the symmetric powers of $\Omega_{Y/S}$  
should be the same as ours and should thus not depend on the dimension of the total space).

One interesting aspect of our result is thus that it removes this indeterminacy. 
However, the main interest of the present text  is the method of proof, which 
is elementary (whereas Franke's and Eriksson's approaches require 
a vast categorical apparatus and use higher $K$-theory, resp. the homotopy theory of schemes). Our canonical isomorphism 
is constructed very explicitly, making it in principle possible to compute 
its norm, when both sides are endowed with metrics (eg Quillen metrics). We hope to return to this question in a later article. 

Note that other constructions of the higher dimensional Deligne pairing were given 
in \cite{Zhang-Heights} and \cite{Elkik-Fibres} but they cannot be used in our context, because they are based on 
cycles classes rather than line bundles and therefore cannot easily be compared 
with our error term. In \cite{BSW-Deligne}, a canonical isomorphism between Ducrot's pairing and Zhang's pairing 
is announced (in a restricted setting), which could be used to bypass the use of Ducrot's pairing in certain situations.

Finally, note that in the situation where $d=1$, Deligne also constructed canonical isomorphisms of line 
bundles (see \cite{Deligne-Det}). Deligne's work was in fact the initial motivation for the work of 
Franke and Eriksson. Under the assumptions of Theorem \ref{mainth} and when $d=1$, Deligne's theorem 
\cite[Th. 9.9 (3)]{Deligne-Det} provides the isomorphism
\begin{equation}
\lambda(\CL)^{\otimes 18}\simeq\lambda(\CO_Y)^{18}\otimes\lambda(\CL^{\otimes 2}\otimes\Omega_{Y/S}^{\vee})^{\otimes 6}\otimes\lambda(\CL\otimes\Omega_{Y/S}^\vee)^{\otimes(-6)}
\label{delq}
\end{equation}
which can be seen as a variant of isomorphism \ref{mainis} when $d=1$. It is not clear to the author whether 
Deligne's theorem implies the existence of isomorphism \ref{mainis} when $d=1$. This is partly a combinatorial problem.

The structure of the article is as follows. In section \ref{secquot} we recall various facts about quotients 
of schemes by finite groups and we prove various supplementary properties of these in the situation 
where the group is a cyclic group of prime order and the fixed point scheme is a Cartier divisor. In section \ref{secducrot} 
we recall the part of Ducrot's work that is relevant to this text. In section \ref{secinv}, we give a proof 
of a local refinement of the fixed formula for an involution, in the situation where the fixed scheme is regularly immersed. 
In the final section \ref{secarr}, we apply this formula to the fibre product of a relative scheme by itself and we prove  
Theorem \ref{mainth}. Note that the core of the proof of Theorem \ref{mainth} amounts to a detailed analysis of the geometry 
of the blow-up of this fibre product along the diagonal. This is intriguing, since this particular space was believed 
to be relevant to a possible solution of the standard conjectures in the early days of scheme theory. It would 
be interesting to relate our construction to statements about algebraic cycles.

{\bf Acknowledgments.} We are grateful to Jean-Michel Bismut and Vincent Maillot for interesting discussions around this article. 

\section{The geometry of quotients by finite groups}
\label{secquot}

 Let $G$ be a finite group. 
 
 A scheme 
$T$ together with a group homomorphism $G\to\Aut(T)$ will be called 
a $G$-equivariant scheme, or an equivariant scheme for short (if there is no 
ambiguity).  A $G$-equivariant morphism of $G$-equivariant schemes 
is a morphism commuting with the action of $G$ on source and target. 
We shall say that the action of $G$ on the $G$-equivariant scheme $T$ is trivial 
if the image of $G\to\Aut(T)$ is the identity morphism. 
 
 A  
$G$-equivariant sheaf (or equivariant sheaf for short) $F$ on a $G$-equivariant scheme is a quasi-coherent sheaf $F$ together with a  
morphism of sheaves $\alpha_g=\alpha_{F,g}:F\to g_*(F)$ for every $g\in G$, such $g_*(\alpha_h)=\alpha_{g\circ h}$ for any $g,h\in G$ 
and $\alpha_{\Id_G}=\Id_F$. 

Suppose that $T$ is a $G$-equivariant scheme with trivial action and that $F$ is a $G$-equivariant sheaf on $T$. The $G$-equivariant structure on $F$ then amounts to 
a homomorphism of groups $G\to\Aut(F)$. We then write $F^G$ for the quasi-coherent sheaf on $T$ such that 
$$
F^G(U)=F(U)^G
$$
for every open set $U\subseteq T$. Here $F(U)^G$ is the subgroup of elements of $F(U)$, which are fixed under the action of $G$.

Suppose that $\phi:T\to Z$ is a morphism of schemes, where $T$ is noetherian. Suppose that $T$ carries 
$G$-equivariant structure and that 
$\phi\circ g=\phi$ for all $g\in G$. Suppose that $F$ is a $G$-equivariant quasi-coherent sheaf. 
Then the sheaf $\phi_*(F)$ is also quasi-coherent. Furthermore, if $Z$ is viewed as a $G$-equivariant scheme carrying the trivial $G$-equivariant structure, then $\phi_*(F)$ carries the $G$-equivariant structure given 
for any $g\in G$ by the composition of 
arrows
$$
\phi_*(F)\stackrel{\to}{\simeq}\phi_*(g_*(F))\stackrel{\to}{\simeq}\phi_*(F)
$$
arising from the equivariant structure on $F$ and the identity $\phi\circ g=\phi.$

Suppose that $\phi:T\to Z$ is a morphism of schemes. Suppose that $T$ carries a $G$-equivariant structure and that 
$\phi\circ g=\phi$ for all $g\in G$. View $Z$ as a $G$-equivariant scheme endowed with the trivial 
$G$-equivariant structure. Let $F$ be a $G$-equivariant quasi-coherent sheaf on $Z$. Then the quasi-coherent sheaf 
$\phi^*(F)$ carries a natural $G$-equivariant structure, given for any $g\in G$ by the composition 
of arrows
$$
\phi^*(F)\stackrel{\to,\phi^*(g_*)}{\simeq}\phi^*(F)\stackrel{\to}{\simeq}g^{-1,*}(\phi^*(F))=g_*(\phi^*(F))
$$
where the first arrow comes by functoriality from the arrow $g_*(F)\to g_*(F)$, the second arrow from the identity $\phi\circ g=\phi$ and the third arrow from the identity $g^{-1,*}=g_*$.

If $x\in X$, then we define $G_d(x)$ to be the stabiliser in $G$ of $x$ viewed as a subset of $X$. This group 
is called the decomposition group of $x$. 
The group $G_d(x)$ naturally acts on the residue field $\kappa(x)$ of $x$. The 
kernel of the homomorphism $G_d(x)\to\Aut(\kappa(x))$ is called the inertia group 
$G_i(x)$ of $x$. 

Suppose that $X$ is a $G$-equivariant scheme. 
A (categorical) quotient $X/G$ of $X$ by $G$ (if it exists) is a $G$-equivariant scheme $X/G$ 
together with an $G$-equivariant morphism \mbox{$q:X\to X/G$,} with the following properties:

- $X/G$ carries the trivial action;

- if $X'$ is a scheme with a trivial $G$-action and 
$q':X\to X'$ is a morphism then there is a unique morphism 
$h:X/G\to X'$, such that $h\circ q=q'$.

These properties clearly determine $X/G$ up to unique isomorphism.

We recall the following 

\begin{prop} Let $X$ be a $G$-equivariant scheme. 
Suppose that the orbit of every point in $X$ is contained in an affine open subscheme. Then the quotient $X/G$ of $X$ by $G$ exists 
and 

\begin{itemize}
\item[\rm (1)] The canonical morphism $q:X\to X/G$ is integral and surjective. 
\item[\rm (2)] The natural morphism of sheaves $\CO_{X/G}\to q_*(\CO_X)$ 
factors through $(q_*(\CO_X))^G$ and induces an isomorphism 
$\CO_{X/G}\to (q_*(\CO_X))^G$. 
\item[\rm (3)] The underlying set of 
$X/G$ is the quotient of the set $X$ by the action of $G$ and the topology 
of $X/G$ is the quotient topology. 
\item[\rm (4)] if 
$Z\to X/G$ is a flat morphism then the natural morphism 
$(Z\times_{X/G}X)/G\to Z$ is an isomorphism.
\item[(5)] Consider the $X/G$-morphism $\phi:G\times_{X/G}X\to X\times_{X/G}X$ 
given in set-theoretic notation by the formula $(g,x)\mapsto(g(x),x)$. Suppose that 
$\phi$ is an isomorphism. Then

- $q$ is \'etale;

- if $\CM$ is  a $G$-equivariant locally free sheaf of finite rank on $X$ then the natural morphism $q^*(q_*\CM)^G\to\CM$ is an isomorphism.
\item[(6)] If $G_i(x)=0$ then $\CO_{X,x}$ is \'etale over $\CO_{X/G,q(x)}$.
\end{itemize}
\label{SGAprop}
\end{prop}
\beginProof See \cite[chap. V, §1 and §2]{SGA1}. \endProof

\begin{cor}
Suppose that there is a morphism of finite type $f:X\to S$, where 
$S$ is a noetherian scheme. Suppose that the action of $G$ on 
$S$ factors through $\Aut_S(X)$. Suppose that the orbit of every point in $X$ is contained in an affine open subscheme.  Then the quotient $X/G$ of $X$ by $G$ exists 
and the morphism $q:X\to X/G$ is finite and surjective.
\label{SGAcor}
\end{cor}

Suppose again that $X$ is a $G$-equivariant scheme. Suppose given a morphism $X\to S$. Suppose that the action of $G$ on 
$X$ factors through $\Aut_S(X)$. In that case, we shall say that 
$X$ is a $G$-equivariant $S$-scheme.

The fixed scheme $X_G$ (if it exists) is a closed subscheme of $X$, which represents the 
functor on $S$-schemes 
$$
T\mapsto X(T)^G.
$$
Note the following link with decomposition and inertia groups: if $x\in X$ and $$G_d(x)=G_i(x)=G$$ then 
$x\in X_G$. This simply follows from the fact that the morphism 
\mbox{$\Spec\,\kappa(x)\to X$} then lies in $X(\Spec\,\kappa(x))^G$.

\begin{prop}
Suppose that $X$ is separated over $S$. Then $X_G$ exists.
\label{corfixex}
\end{prop}
\beginProof Let $\Gamma_g$ be the graph of $g\in G$ in $X\times_S X$. From the separatedness assumption, each $\Gamma_g$ is a closed subscheme of $X\times_S X$. It follows from the definitions that 
we can take $X_G=\cap_{g\in G}\Gamma_g$.
\endProof

If $X_G$ exists, we shall write $N_{X_G/X}$ for the conormal sheaf of $X_G$ in $X$. 
Recall that If $\CI$ is the ideal sheaf of $X_G$ in $X$, we have by definition $N_{X_G/X}=
\CI/\CI^2$. The sheaf $\CI/\CI^2$ has a natural structure of $\CO_{X_G}$-module. 
The conormal sheaf $N_{X_G/X}$ is thus a quasi-coherent sheaf on $X_G$ 
and it carries a natural action of $G$. 

Suppose now that $X$ is a $G$-equivariant $S$-scheme. Suppose that $G_S\simeq\mn_S$, 
where \mbox{$\mn=\Spec\,\mZ[t]/(1-t^n)$} is the diagonalisable group scheme associated with 
the cyclic group $\mZ/n\mZ$. Note that the condition $G_S\simeq\mn_S$ is equivalent to requiring $n$ to be invertible 
in $S$ and for the polynomial $x^n-1$ to split into linear factors in $\Gamma(S,\CO_S).$ 
We note the following two facts. 

Suppose in this paragraph that $X=\Spec\,R$ is affine. Then the action of $G$ on $X$ 
is given by a ring grading $R\simeq\oplus_{k\in\mZ/n\mZ}R_k$, such that 
the morphism $X\to S$ factors through $\Spec\,R_0$. Furthermore, the ideal of $X_G$ is then $R\cdot R_{\not=0}$, where 
$$
R_{\not=0}:=\oplus_{k\in\mZ/n\mZ,\,k\not=0}R_k.
$$
See \cite[proof of Prop. 3.1]{Thomason-Formule} (this is also a good exercise for the reader). 

Suppose that the action of $G$ on $X$ is trivial. Let $F$ be a $G$-equivariant sheaf on $X$. The $G$-equivariant structure on 
$F$ is then given by a $\mZ/n\mZ$-grading of $\CO_X$-modules $$F\simeq\oplus_{k\in\mZ/n\mZ} F_k.$$ The action of 
$G=\mn(S)$ on $F$ is then by construction given  by the formula 
$$\rho(\oplus_{k\in\mZ/n\mZ}\,f_k)=\oplus_{k\in\mZ/n\mZ}\,\rho^k\cdot f_k,$$
where $\rho\in \mn(S)$ and $f_k$ is a local section of $F_k$. In particular, we have 
$F_0=F^G$. 

We also record the following

\begin{lemma}
Let $X$ be an $G$-equivariant $S$-scheme. Suppose that the orbit of every point in $X$ is contained in an affine open subscheme. Suppose that $G_S\simeq\mn_S$. 
If 
$Z\to X/G$ is a morphism then the natural morphism 
$(Z\times_{X/G}X)/G\to Z$ is an isomorphism.
\label{univbc}
\end{lemma}
In other words, when $G_S\simeq\mn_S$, the quotient construction commutes with any base change 
on $X/G$ (not only flat base changes as in Proposition \ref{SGAprop} (4)).
\beginProof
By Proposition \ref{SGAprop} (4), we may assume that $Z$ and $X$ are affine, say 
$Z=\Spec\,B$ and $X=\Spec\,A$. In this case, we have to prove that 
the morphism of $A_0$-modules
$$
B\to (B\otimes_{A_0}A)_0
$$
given by the formula $b\mapsto b\otimes 1$ 
is an isomorphism. We have 
$$
B\otimes_{A_0}A=B\otimes_{A_0}\bigoplus_{k\in\mZ/n\mZ}A_k=\bigoplus_{k\in\mZ/n\mZ}B\otimes_{A_0}A_k
$$
so that $(B\otimes_{A_0}A)_0=B\otimes_{A_0}A_0=B$, proving the assertion.\endProof

\begin{prop}  Suppose that $X$ is a $G$-equivariant $S$-scheme such that $S$ is noetherian and the morphism 
$X\to S$ is separated and of finite type. Suppose that the orbit of every point in $X$ is contained in an affine open subscheme. Finally, suppose 
that $G_S\simeq\mn_S$. Let $\iota:X_G\to X$ be the fixed point scheme of $X$. 
 Then:
\begin{itemize}
\item[\rm (1)] Suppose that $n$ is prime and that $X_G$ is a  (possibly empty) Cartier divisor. Then $q$ is flat.
\item[\rm (2)] Suppose that $X_G$ is a  Cartier divisor. Then $(N_{X_G/X})_0=0$.
\item[\rm (3)] The morphism $q\circ\iota:X_G\to X/G$ is a closed immersion and we have the set-theoretic equality 
$q^{-1}(q(X_G))=X_G$. Thus we have a natural isomorphism $(X/G)\without q(X_G)\simeq
(X\without X_G)/G$.
\item[\rm (4)] Let $U=X\without X_G$ (so that $U/G=(X/G)\without q(X_G)$ by (3)). Consider the $U/G$-morphism $$\phi:G\times_{U/G}U\to U\times_{U/G}U$$ 
given in set-theoretic notation by the formula $(g,u)\mapsto(g(u),u)$. If $n$ is prime then
$\phi$ is an isomorphism.
\item[(5)] Let $\CM$ be a $G$-equivariant locally free sheaf of finite rank on $X$. 
Suppose that $\iota^*\CM$ carries the trivial action, that $q$ is flat and that $n$ is  prime. Then the natural morphism
$q^*(q_*\CM)_0\to\CM$ is an isomorphism.
\item[(6)] If $X\to S$ is smooth and $X_G\to S$ is flat then $X_G\to S$ is smooth.
\item[(7)] If $X\to S$ is smooth, $X_G$ is a Cartier divisor in $X$ and $X_G\to S$ is 
flat then  
$X/G\to S$ is also smooth.
\end{itemize}
\label{vimplem}
\end{prop}
\begin{rem}\rm A variant (for algebraic varieties) of (5) is proven in \cite[Th. 2.3]{DN-Groupe}.\end{rem}
\beginProof We begin with (1). We may suppose that $X=\Spec(R)$ is affine. 
Then $X/G=\Spec(R_0)$ by Proposition \ref{SGAprop} (2). To show that $R$ is flat over $R_0$, it is sufficient 
to show that for all $\mfp\in\Spec(R)$, the ring $R_\mfp$ is flat 
over the ring $R_{0,\mfp\cap R_0}$. If $\mfp\not\supseteq R\cdot R_{\not=0}$, then 
$\mfp\not\in X_G$ by the previous discussion. Thus $G_i(x)\not=G$ and thus 
$G_i(x)=0$ since $n$ is prime; thus $R_\mfp$ is flat 
over the ring $R_{0,\mfp\cap R_0}$ by Proposition \ref{SGAprop} (6). Thus we may assume that $\mfp\supseteq R\cdot R_{\not=0}$. The prime ideal 
$\mfp$ is then graded by construction (if $r\in \mfp$, write 
$r=r_0+\dots+r_{n-1}$, where the $r_i$ are homogenous for the grading; by assumption  $r_1,\dots,r_{n-1}\in\mfp$; thus 
$r_0\in\mfp$ as well). 
The ring $R_\mfp$ is thus naturally a $\mZ/n\mZ$-graded local ring. 
Now notice that we have a natural identification
$$
R_{0,\mfp\cap R_0}=(R_\mfp)_0
$$
(use the fact that $R\backslash\mfp\subseteq R_0$). Also by construction the ideal generated by the image of the ideal  $R\cdot R_{\not=0}$ in 
$R_\mfp$ is $R_\mfp\cdot R_{\mfp,\not=0}$. Thus the assumption that $R\cdot R_{\not=0}$ is a Cartier divisor implies that there exists 
$t\in R_\mfp$, which is not a zero divisor, such that $(t)=R_\mfp\cdot R_{\mfp,\not=0}$.

Thus we may assume without restriction of generality that 
$R$ is a local ring and that $R\cdot R_{\not=0}$ is generated 
by an element $t$, which is not a zero divisor.  

We claim that $t$ can be taken to be homogenous of 
degree $\not=0\ ({\rm mod}\ n)$.

To verify the claim, let 
$$
R\cdot R_{\not=0}=(a_1,\dots, a_k)
$$
where the $a_i\in R_{\not=0}$ are homogenous (recall that $R$ is noetherian). 
We take $k$ minimal. We may assume that $k>1$, otherwise there is nothing to prove. 
Then for some family of $x_i\not=0$, we have
$$
x_1 a_1+\cdots+x_k a_k=t
$$
Let $b_1\in R$ be such that $a_1=t\cdot b_1$. If $b_1$ is a unit then 
$R\cdot R_{\not=0}=(a_1)$ contradicting the assumption that $k>1$. 
Thus $b_1$ is not a unit and thus $1-x_1 b_1$ is a unit since $R$ is local. We compute
$$
t={a_2 x_2\over 1-x_1b_1}+\dots +{a_k x_k\over 1-x_1b_1}
$$
contradicting minimality again. Thus $k=1$ and the claim is verified.

So we may suppose that $(t)=R\cdot R_{\not=0}$ where $t$ is homogenous of degree $\not=0\ ({\rm mod}\ n)$. 
Now note that $t^n\in R_0$. By the local criterion of flatness (see \cite[chap. 8, §22]{Matsumura-Commutative}), to verify that 
$R$ is flat over $R_0$, it is sufficient to verify that $t^n$ is not a zero-divisor 
in $R_0$ and that $R/(t^n)$ is flat over $R_0/(t^n)$. The  first condition is satisfied by construction. To verify that $R/(t^n)$ is flat over $R_0/(t^n)$ 
note that $R/(t^n)$ has a the finite filtration 
\begin{equation}
R/(t^n)\supseteq (t)/(t^n)\supseteq(t^2)/(t^n)\supseteq\dots\supseteq(t^{n-1})/(t^n)\supseteq 0
\label{filteq}
\end{equation} whose quotients are isomorphic to $R/(t)$ (the fact that $t$ is not 
a zero divisor is used here). It is thus sufficient to show that $R/(t)$ is 
flat over $R_0/(t^n)$ via the natural map $R_0/(t^n)\to R/(t)$. For this, 
note that since $t^n$ is of degree $0\ ({\rm mod}\ n)$ we have a natural isomorphism
$$
R_0/(t^n)\simeq (R/(t^n))_0.
$$
Furthermore the degree
of $t^i$ in $R$ for $i=0,1\dots, n-1\ ({\rm mod}\ n)$ is $i\cdot\deg(t)\ ({\rm mod}\ n)$ and thus the degrees 
of the $t^i$ for $i=0,1\dots, n-1$ are all distinct, since $n$ is prime. Hence 
the filtration \refeq{filteq} splits and we have 
we have an isomorphism of graded rings $$R/(t^n)\simeq\oplus_{i=0}^{n-1}(t^i)/(t^n).$$
In particular $(R/(t^n))_0\simeq R/(t)$ and thus $R/(t)$ is 
flat over $R_0/(t^n)\simeq(R/(t^n))_0\simeq R/(t).$

To prove (2), localising at points of $X_G$, we may still assume that $X=\Spec(R)$, where $R$ is a local ring 
and $R\cdot R_{\not=0}$ is generated by a single element $t$, which is not a zero divisor.
In the proof of (1), it was shown that we may suppose that $t$ is homogenous of degree $\not=0$. The sheaf  
$N_{X_G/X}$ corresponds to the $R$-modules $(t)/(t^2)$ and thus 
$(N_{X_G/X})_0=0$, since $t$ is of degree $\not=0\ ({\rm mod}\ n)$.

Proof of (3). We may suppose that $X=\Spec\,R$ is affine. The first statement 
now corresponds to the statement that $R_0\to R/(R\cdot R_{\not=0})$ is surjective. 
This follows from the definitions. The fact that $q^{-1}(q(X_G))=X_G$ follows from 
Proposition \ref{SGAprop} (3). The third assertion follows from Proposition \ref{SGAprop} (4).

Proof of (4). Note that for all $x\in X\without X_G$, we have $G_i(x)\not=G$ and 
thus $G_i(x)=0$, since $n$ is prime.  By Proposition \ref{SGAprop} (6) this implies that $q$ is \'etale, in particular flat. Hence the morphism $U\to
U/G$ is finite and flat. 

We first compute its degree. For this, let $u_0\in U/G$ and 
let $H$  be the spectrum of the strict henselisation of $\CO_{U/G,u_0}$. Then $H\simeq(U\times_{U/G}H)/G$ by  Proposition \ref{SGAprop} (4) and the 
fact that $H$ is flat over $\CO_{U/G,u_0}$ (see \cite[I, §1, 1.20]{FK-Etale} for this). We only have to compute the degree of $U\times_{U/G}H$ over $H$. Now note that $U\times_{U/G}H$ 
is a disjoint union $\coprod_{i\in I}H_i$ of copies of $H$, since $H$ is strictly henselian and 
$U\times_{U/G}H\to H$ is \'etale. Furthermore, the group $G$ permutes the $H_i$ 
and also the closed points of the $H_i$. Hence the degree is the cardinality of the 
orbit of a closed point $P\in H_{i_0}$ ($i_0$ arbitrary). Since $G_i(P)=G_d(P)$, we must have $G_d(P)=0$, since $n$ is prime and $(U\times_{U/G}H)_G$ is empty. Hence the orbit of 
$P$ has $n$ elements and thus the degree of $U\to
U/G$ is $n$.

Now consider the morphism $\phi:G\times_{U/G}U\to U\times_{U/G}U$. Let $T$ be a connected scheme. The map 
$G(T)\times_{(U/G)(T)} U(T)\to U(T)\times_{(U/G)(T)}U(T)$ is injective. To see this note that otherwise 
there is $e\in U(T)$ and $g\in G(T)$ such that $g\not=0$ and $g(e)=e$; since 
$G(T)$ is of prime order this means that $e\in U(T)^G$ and thus $e\in U_G(T)$, which 
is not possible, since $U_G$ is empty. Since $T$ was arbitrary, the morphism $\phi$ is a monomorphism 
of schemes. Since it is also proper (because $G\times_{U/G}U$ and $U\times_{U/G}U$ are proper over $U/G$), it is a closed immersion (see \cite[IV.3, 8.11.5]{EGA} for this). Since both 
$G\times_{U/G}U$ and $U\times_{U/G}U$ are flat and finite of the same rank over $U$ by the previous paragraph, this implies that 
$\phi$ is an isomorphism.

Proof of (5). Consider the natural morphism
$$
\alpha:q^*(q_*\CM)_0\to\CM
$$
The restriction $\alpha_{X\backslash X_G}$ is an isomorphism 
by (4) and Proposition \ref{SGAprop} (5). Since both sides are locally free of finite rank, by Nakayama's lemma, it is sufficient to 
show that $\alpha_{\kappa(x)}$ is surjective for $x\in X_G$. In particular, it is sufficient 
to show that the restriction $\iota^*(\alpha)$ of $\alpha$ to $X_G$ is an isomorphism. 
Now note that since $q$ is an affine morphism, the natural adjunction morphism $
\alpha:q^*(q_*\CM)\to\CM$ is a surjection and thus we have a surjection
$$
\iota^*(q^*(q_*\CM))\to\iota^*(\CM)
$$
extending $\iota^*(\alpha).$ Hence we have a surjection
$$
\iota^*(q^*((q_*\CM)_0))\to\iota^*(\CM)_0
$$
and since $\iota^*(\CM)_0=\iota^*(\CM)$ by assumption we get a surjection
$$
\iota^*(q^*((q_*\CM)_0))\simeq\iota^*(\CM)
$$
which must be an isomorphism, since both sides are locally free of the same rank.

Proof of (6). We need to check that the geometric 
fibres $X/G\to S$ are regular. So let $\Spec\, k\to S$ 
be a geometric point. By assumption, $X_k$ is regular and by 
\cite[Prop. 3.1]{Thomason-Formule} , $(X_k)_G=(X_G)_k$ is then also regular.

Proof of (7). Since $q$ is faithfully flat, we see that $X/G\to S$ is also flat.
To see that $X/G\to S$ is smooth, we need to check that the geometric 
fibres $X/G\to S$ are regular. Now since $X_G$ is flat over $S$ and a Cartier divisor, we see that for 
any base change $T\to S$, $(X_T)_G\to T$ is also flat and a Cartier divisor. 
Furthermore, by Lemma \ref{univbc}, for 
any base change $T\to S$, we have $(X/G)_T\simeq (X_T)/G$. So let $\Spec\, k\to S$ 
be a geometric point. By assumption $X_k$ is regular and since $(X_k)_G$ is 
a Cartier divisor, we see that $(X_k)/G=(X/G)_k$ is regular, since $q_k$ is faithfully flat by (1) 
and Proposition \ref{SGAprop} (1). 
\endProof

\section{Ducrot's generalisation of the Deligne pairing}
\label{secducrot}

Let $g:Y\to S$ be a smooth and strongly projective morphism of constant relative dimension $d$. Suppose that $S$ is noetherian. If $F_1,\dots,F_k$ is a finite sequence of coherent locally free sheaves on $X$, we shall write 
$$
\lambda(n_1 F_1+\dots +n_k F_k):=
\bigotimes_{r=1}^k\lambda(F_r)^{\otimes n_r}
$$
for any $n_1,\dots,n_k\in\mZ.$

Let $\CL_1,\dots\CL_{d+1}$ be line bundles on $Y$. Ducrot showed in 
\cite[§5]{Ducrot-Cube} that the line bundle 
$$
I_{Y/S}(\CL_1,\dots,\CL_{d+1}):=\lambda((\CO_Y-\CL_1)\otimes(\CO_Y-\CL_2)\otimes\dots\otimes(\CO_Y-\CL_{d+1}))
$$
is multiadditive in the line bundles $\CL_1,\dots\CL_{d+1}$. In particular, he shows 
that if $\CQ$ is a line bundle on $Y$, then there is a canonical isomorphism
\begin{equation}
I_{Y/S}(\CL_1\otimes\CQ,\dots,\CL_{d+1})\simeq
I_{Y/S}(\CL_1,\dots,\CL_{d+1})\otimes
I_{Y/S}(\CQ,\dots,\CL_{d+1})
\label{multadd}
\end{equation}
It is very plausible that the canonical isomorphism \refeq{multadd} is compatible with any base change to a noetherian scheme. However, we were not able to verify this in detail. 

We shall call {\it hypothesis} (H) the statement that the canonical isomorphism \refeq{multadd} is compatible with any base change to a noetherian scheme. 

We may thus compute
\begin{eqnarray*}
&&\lambda((\CO_Y-\CQ)\otimes(\CO_Y-\CL_1)\otimes(\CO_Y-\CL_2)\otimes\dots\otimes(\CO_Y-\CL_{d+1}))
\\&\simeq&
\lambda((\CO_Y-\CL_1)\otimes(\CO_Y-\CL_2)\otimes\dots\otimes(\CO_Y-\CL_{d+1}))
\\&\otimes&\lambda((\CQ-\CQ\otimes\CL_1)\otimes(\CO_Y-\CL_2)\otimes(\CO_Y-\CL_3)\otimes\dots\otimes(\CO_Y-\CL_{d+1}))^{\vee}
\\&\simeq&
\lambda((\CO_Y-\CL_1)\otimes(\CO_Y-\CL_2)\otimes\dots\otimes(\CO_Y-\CL_{d+1}))
\\&\otimes&\lambda(\big(\CO_Y-\CL_1\otimes\CQ-(\CO_Y-\CQ)\big)\otimes(\CO_Y-\CL_2)\otimes(\CO_Y-\CL_3)\otimes\dots\otimes(\CO_Y-\CL_{d+1}))^{\vee}
\\&\simeq&
I_{Y/S}(\CL_1),\dots,\CL_{d+1}))
\\&\otimes&
I_{Y/S}(\CL_1\otimes\CQ,\CL_2,\dots,\CL_{d+1})^\vee\otimes
I_{Y/S}(\CQ,\CL_2,\dots,\CL_{d+1})
\\&\simeq&
I_{Y/S}(\CL_1,\CL_2,\dots,\CL_{d+1})\otimes
I_{Y/S}(\CL_1,\CL_2,\dots,\CL_{d+1})^\vee
\\&\otimes&I_{Y/S}(\CQ,\CL_2,\dots,\CL_{d+1})^\vee\otimes
I_{Y/S}(\CQ,\CL_2,\dots,\CL_{d+1})\simeq\CO_X
\end{eqnarray*}
and conditional on hypothesis (H) this trivialisation is invariant under any base-change to a noetherian scheme.

The following proposition summarises the discussion.

\begin{theor}
Suppose that $Y\to S$ is smooth and strongly projective. Suppose that $S$ is noetherian. Let 
$\CL_1,\dots,\CL_{d+2}$ be line bundles on $Y$. Then the line 
bundle 
$$
\lambda((\CO_Y-\CL_1)\otimes(\CO_Y-\CL_2)\otimes\dots\otimes(\CO_Y-\CL_{d+2}))
$$
is canonically trivial and conditional on hypothesis (H) this trivialisation is invariant under base-change to a noetherian scheme.
\label{DucrotT}
\end{theor}

\section{Local refinement of the fixed point formula for an involution}
\label{secinv}

Let $S$ be a noetherian scheme and let $f:X\to S$ be a flat and separated morphism 
of finite type. Suppose that $2$ is invertible in $S$. 
Let $G=\mZ/2$, so that we have canonical isomorphism 
$G_S\simeq\m2_S$. Suppose that  we have a $G$-equivariant structure 
on $X$ over $S$. Suppose finally that the orbit of every point in 
$X$ is contained in an open affine subscheme. Let \mbox{$\iota:X_G\hookrightarrow X$} be the fixed scheme of $X$ and let  
$q:X\to X/G$ be the quotient morphism. These morphisms exist by Proposition \ref{corfixex} and 
Theorem \ref{SGAprop}. 
Note that if $q$ is flat then it is faithfully flat (since it is surjective) and thus if 
$q$ is flat the natural morphism $X/G\to S$ is also 
flat. Similarly, if $f$ is strongly projective then so is the natural morphism 
$X/G\to S$. 

If $F$ is a quasi-coherent sheaf on $X$, we shall say that $F$ is $f$-acyclic is 
$\R^k f_*(F)=0$ when $k>0$. If the action on $X$ is trivial and $F$ is an equivariant locally free sheaf on 
$X$, we shall often write $F_0=F_+$ and $F_-=F_1$.

From now on, we suppose that $f$ is strongly projective as well. 

\begin{lemma}  
If $J^\bullet$ is a bounded complex of $G$-equivariant coherent sheaves on 
$X$, then there is a bounded complex $H^\bullet$ of $G$-equivariant $f$-acyclic coherent sheaves on $X$ and a $G$-equivariant quasi-isomorphism 
$J^\bullet\to H^\bullet$. If the sheaves $J^i$ are locally free then 
$H^\bullet$ can be chosen so that the sheaves $H^i$ are also locally free.
\label{lemglob}
\end{lemma}
\beginProof Let $F$ be a $G$-equivariant coherent sheaf on $X$. Let $\CM$ be a relatively ample line bundle on $X$. This exists because $f$ is strongly projective over $S$. 
Let $g_0$ be the unique generator of $G$. The line bundle $\CA:=\CM\otimes g_{0,*}(\CM)$ is then 
naturally $G$-equivariant. The line bundle $\CA$ is also ample and thus there is a natural number $n_0>0$ such that 
$$\R^k f_*(F\otimes\CA^{\otimes n})=0\,\,\,\,{\rm and}\,\,\,\,\R^k f_*(\CA^{\otimes n})=0 $$ for all $k>0$ and for all $n\geq n_0$ and such that the natural morphism 
$$f^*(f_*(\CA^{\otimes n_0}))\mapsto \CA^{\otimes n_0}$$ is surjective. 
Note that by the semicontinuity theorem, $f_*(\CA^{\otimes n_0})$ is then a locally free sheaf.
 
Let 
$r:=\rk(f^*(f_*(\CA^{\otimes n_0})))$. The previous morphism induces an 
exact $G$-equivariant Koszul resolution
$$
\hskip-1cm
0\to \Lambda^{r}(f^*(f_*(\CA^{\otimes n_0})))\otimes\CA^{\otimes(-rn_0)}\to\dots\to 
\Lambda^{2}(f^*(f_*(\CA^{\otimes n_0})))\otimes\CA^{\otimes(-2n_0)}\to f^*(f_*(\CA^{\otimes n_0}))\otimes\CA^{\otimes(-n_0)}\to\CO_X\to 0.
$$
Dualising this complex and tensoring by $F$, we get the exact $G$-equivariant complex
$$
0\to F\to f^*(f_*(\CA^{\otimes n_0})^\vee)\otimes\CA^{\otimes n_0}\otimes F\to 
\Lambda^2(f^*(f_*(\CA^{\otimes n_0})^\vee))\otimes\CA^{\otimes 2n_0}\otimes F\to\dots.
$$
We have thus constructed a finite $G$-equivariant resolution of $F$ by $f$-acyclic coherent  sheaves. If $F$ is locally free, the resolution will also consist of locally free sheaves. This proves the lemma in the situation where $J^\bullet$ consists of 
one object concentrated in degree $0$. Now suppose that $F'$ is another coherent $G$-equivariant 
sheaf on $X$ and that $F\to F'$ is a $G$-equivariant morphism of sheaves. 
We may repeat the above construction for $F'$, choosing an $n_0$ which is sufficiently large so that it can used for both $F$ and $F'$. One then obtains a commutative diagram with exact rows
\begin{diagram}
0&\rTo& F&\rTo& f^*(f_*(\CA^{\otimes n_0})^\vee)\otimes\CA^{\otimes n_0}\otimes F&\rTo& \Lambda^2(f^*(f_*(\CA^{\otimes n_0})^\vee))\otimes\CA^{\otimes 2n_0}\otimes F&\rTo&\dots\\
 & &\dTo& &\dTo& &\dTo\\
0&\rTo& F'&\rTo& f^*(f_*(\CA^{\otimes n_0})^\vee)\otimes\CA^{\otimes n_0}\otimes F'&\rTo& 
\Lambda^2(f^*(f_*(\CA^{\otimes n_0})^\vee))\otimes\CA^{\otimes 2n_0}\otimes F'&\rTo&\dots
\end{diagram}
Generalising this to a complex of $G$-equivariant sheaves, we may associate a double complex 
of $f$-acyclic coherent  sheaves with $J^\bullet$. The total complex of this double complex will be quasi-isomorphic to $J^\bullet$ and it will consists of locally free sheaves if 
$J^\bullet$ consists of locally free sheaves. We leave the details to the reader.
\endProof
Let now $\Coh^\eq(X)$ (resp. $\Coh^\eq(S))$) be the category of coherent 
$G$-equivariant sheaves on $X$ (resp. on $S$). These categories carry natural structures of abelian categories and the functor $f_*$ induces a left exact functor 
from $\Coh^\eq(X)$ to $\Coh^\eq(S)$, that we shall 
call $f^\eq_*$ to underline the dependence on the equivariant structures of 
$X$ and $S$. In view of Lemma \ref{lemglob} 
and \cite[Th. I.5.1]{Hartshorne-Residues}, the functor $f^\eq_*$ has 
a right derived functor $$\R^\bullet f^\eq_*:D^b(\Coh^\eq(X))\to D^b(\Coh^\eq(S)).$$
where $D^b(\Coh^\eq(X))$ (resp. $D^b(\Coh^\eq(S))$) is the derived category of bounded complexes 
of $G$-equivariant coherent sheaves on $X$ (resp. on $S$).

Denote by $\Coh(X)$ (resp. $\Coh(S))$) the category of coherent 
sheaves on $X$ (resp. on $S$). The functors $(\bullet)_{-}=(\bullet)_{1}:\Coh^\eq(S)\to \Coh(S)$ and 
$(\bullet)_{+}=(\bullet)_{0}:\Coh^\eq(S)\to \Coh(S)$ are 
exact functors and so they uniquely extend to functors from 
$D^b(\Coh^\eq(S))$ to $D^b(\Coh(S))$, which are their right and left derived functors simultaneously. We shall also call these extensions $(\bullet)_{+}$ and $(\bullet)_{-}$. 

Let now $F$ be a coherent locally free $G$-equivariant sheaf on $X$. 
By Lemma \ref{lemglob}, the object $\R^\bullet f^\eq_*(F)$ is represented by a bounded complex 
of $G$-equivariant locally free sheaves and thus the objects 
$(\R^\bullet f^\eq_*(F))_{-}$ and $(\R^\bullet f^\eq_*(F))_{+}$ of $D^b(\Coh(S))$ are perfect complexes.

We define
$$
\lambda(F):=\det((\R^\bullet f^\eq_*(F))_{+})\otimes \det((\R^\bullet f^\eq_*(F))_{-})^\vee
$$
where $\det(\bullet)$ is the Knudsen-Mumford determinant of a perfect complex 
(see \cite{KM-Det}). 

 Note that with this definition, if 
\begin{equation}
0\to F'\to F\to F''\to 0
\label{ExF}
\end{equation}
is an exact sequence of $G$-equivariant coherent locally free sheaves, we have a canonical isomorphism 
\begin{equation}
\lambda(F')\otimes\lambda(F'')\simeq\lambda(F). 
\label{detdef}
\end{equation}
Indeed, the sequence \refeq{ExF} defines a triangle in $D^b(\Coh^\eq(X))$ 
and thus induces a triangle
$$
\R^\bullet f^\eq_*(F')\to \R^\bullet f^\eq_*(F)\to \R^\bullet f^\eq_*(F'')\to \R^\bullet f^\eq_*(F')[1]
$$
in $D^b(\Coh^\eq(S))$. Thus we obtain two triangles
$$
\R^\bullet f^\eq_*(F')_\pm\to \R^\bullet f^\eq_*(F)_\pm\to \R^\bullet f^\eq_*(F'')_\pm\to \R^\bullet f^\eq_*(F')_\pm[1]
$$
and we have canonical ismorphisms
$$
\lambda(\R^\bullet f^\eq_*(F)_\pm)\simeq \lambda(\R^\bullet f^\eq_*(F')_\pm)\otimes
\lambda(\R^\bullet f^\eq_*(F'')_\pm)
$$
from the standard properties of the determinant functor. This 
shows that \refeq{detdef} holds. The identity \refeq{detdef} makes sense more 
generally if $F', F$ and $F'$ are $G$-equivariant coherent sheaves, which 
have the property that $\R^\bullet f^\eq_*(F')$, $\R^\bullet f^\eq_*(F)$ and $\R^\bullet f^\eq_*(F'')$ 
can be represented by bounded complexes of $G$-equivariant locally free sheaves. 

If $F_1,\dots,F_k$ is a finite sequence of equivariant coherent locally free sheaves on $X$, we shall write 
$$
\lambda(n_1 F_1+\dots+n_k F_k):=
\bigotimes_{r=1}^k\lambda(F_r)^{\otimes n_r}
$$
for any $n_1,\dots,n_k\in\mZ.$

Finally we shall write $\TW$ for the trivial sheaf $\CO_X$, endowed 
with the $G$-equivariant structure such that for any $\rho\in\m2(S)$ the 
isomorphism $\alpha_{\rho,\TW}:
\TW\to g_*(\TW)$ composed with the canonical non equivariant identification $g_*(\TW)\simeq\TW$ is given by multiplication by $\rho$.
If $F$ is an equivariant sheaf on $X$, we shall write 
$F\TW$ for $F\otimes\TW$. Note that if $F$ is an equivariant coherent locally 
free sheaf on $X$, we have 
\begin{equation}
\lambda(F\{-1\})\simeq\lambda(F)^{\vee}.
\label{veeq}
\end{equation}

In this section, we shall prove a version of the relative geometric fixed point formula for the $G$-action of $G$ on $X$, which avoids $K$-theory entirely, 
replacing all the equalities in a Grothendieck group or a Picard group by explicit isomorphisms:

\begin{theor}
In addition to the above assumptions, suppose that $f$ is smooth.  Let $\CM$ be a $G$-equivariant coherent locally free sheaf of rank one on $X$. Suppose that $f$ has constant relative dimension $d$. Suppose also that the morphism 
$X_G\to S$ is flat. Then $X_G\to S$ is smooth and thus 
$X_G$ is regularly immersed in $X$. Let $N=N_{X_G/X}$ be the conormal 
bundle of $\iota$, endowed with its canonical $G$-equivariant structure. We have a canonical isomorphism
\begin{equation}
\lambda(\CM)^{\otimes 2^{d+1}}\simeq 
 \lambda\Big(\iota^*(\CM)\otimes
 \sum_{i=0}^d\sum_{j=0}^i 2^{d-i}{i\choose j}\Sym^j(N)\Big).
 \end{equation}
Conditional on hypothesis (H), this isomorphism is compatible 
with any base-change $h:S'\to S$ such that 
$S'$ is noetherian. 
\label{invfunc}
\end{theor}
For the proof, we shall need the following
\begin{lemma}
Let $Z\to T$ be a morphism of noetherian schemes. 
Let $C\hookrightarrow Z$ be a regular closed immersion. Suppose that 
$C$ is flat over $T$. Let $h:T'\to T$ be a morphism of schemes, where 
$T'$ is noetherian. Then the natural 
morphism $\Bl_{C_{T'}}(Z_{T'})\to\Bl_{C}(Z)_{T'}$ is an isomorphism.
\label{lemibl}
\end{lemma}
\beginProof
Left to the reader. 
\endProof
\beginProof (of Theorem \ref{invfunc}). Suppose first that $X_G$ is a Cartier divisor. 
Let $\CL:=\CO(-X_G)$. 

We have an exact sequence
\begin{equation}
0\to \CL\otimes\CM\to \CM\to \iota_*(\iota^*(\CM))\to 0
\label{exseq}
\end{equation}


Note that by the adjunction formula we have a canonical equivariant isomorphism \mbox{$\iota^*(\CL)\simeq N$.} 
Let $\CJ:=q_*(\CL\{-1\})_0$. By the adjunction formula and Proposition \ref{vimplem} (5) this is a line bundle on $X/G$ such that 
$q^*(\CJ)=\CL\{-1\}$. 

We first list some  identities in $\mQ(t)$. We have
$$
{1\over t}={1\over 2-(2-t)}={1/2\over 1-(2-t)/2}=
{1\over 2}+{(2-t)\over 2^2}+{(2-t)^2\over 2^3}+\dots+{(2-t)^k\over 2^{k+1}}+
{1\over 2}{((2-t)/2)^{k+1}\over 1-(2-t)/2}
$$
so that 
$$
t\cdot[{1\over 2}+{(2-t)\over 2^2}+{(2-t)^2\over 2^3}+\dots+{(2-t)^k\over 2^{k+1}}]=
1-((2-t)/2)^{k+1}
$$
and
$$
t\cdot[2^k+{2^{k-1}(2-t)}+{2^{k-2}(2-t)^2}+\dots+{(2-t)^k}]=
2^{k+1}-(2-t)^{k+1}
$$
in $\mZ[t]$. Define
$$
P_k(t):=2^k+{2^{k-1}(2-t)}+{2^{k-2}(2-t)^2}+\dots+{(2-t)^k}\in\mZ[t]
$$
so that by the above we have $t\cdot P_k(t)=2^{k+1}-(2-t)^{k+1}$. 

Now we compute 
\begin{eqnarray*}
&&\lambda(\iota^*(\CM)\otimes P_k(\CO_{X_G}-\iota^*(\CL))\stackrel{(a)}{\simeq}
\lambda(\iota^*(\CM)\otimes P_k(\CO_{X_G}+\iota^*(\CL\{-1\}))\\
&\stackrel{(b)}{\simeq}&\lambda(\iota^*(\CM)\otimes P_k(\CO_{X_G}+N\TW))\stackrel{(c)}{\simeq}
\lambda(\CM\otimes(\CO_X-\CL)\otimes P_k(\CO_X-\CL))\\&\stackrel{(d)}{\simeq}&
\lambda(\CM\otimes(\CO_X^{\oplus 2^{k+1}}-(\CO_X^{\oplus 2}-(\CO_X-\CL))^{\otimes(k+1)}))
\\&\stackrel{(e)}{\simeq}&
\lambda(\CM\otimes(\CO_X^{\oplus 2^{k+1}}-(\CO_X^{\oplus 2}-(\CO_X+\CL\{-1\}))^{\otimes(k+1)}))
\\&\stackrel{(f)}{\simeq}&
\lambda(\CM\otimes(\CO_X^{\oplus 2^{k+1}}-(\CO_X-\CL\{-1\})^{\otimes(k+1)}))
\\&\stackrel{(g)}{\simeq}&
\lambda(\CM)^{\otimes2^{k+1}}\otimes\lambda(\CM\otimes(\CO_X-\CL\{-1\})^{\otimes(k+1)})^{\vee}
\\&\stackrel{(h)}{\simeq}&
\lambda(\CM)^{\otimes2^{k+1}}\otimes
\lambda(q_*(\CM)\otimes(\CO_{X/G}-J)^{\otimes(k+1)})^{\vee}
\\&\stackrel{(i)}{\simeq}&
\lambda(\CM)^{\otimes2^{k+1}}\otimes
\lambda((q_*(\CM)_+-q_*(\CM)_-)\otimes(\CO_{X/G}-J)^{\otimes(k+1)})^{\vee}
\\&\stackrel{(j)}{\simeq}&
\lambda(\CM)^{\otimes2^{k+1}}\otimes
\lambda(((\CO_{X/G}-q_*(\CM)_-)-(\CO_{X/G}-q_*(\CM)_+))\otimes(\CO_{X/G}-J)^{\otimes(k+1)})^{\vee}
\\&\stackrel{(k)}{\simeq}&
\lambda(\CM)^{\otimes2^{k+1}}\otimes\lambda((\CO_{X/G}-q_*(\CM)_-)\otimes(\CO_{X/G}-J)^{\otimes(k+1)})^{\vee}\\
&\otimes&\lambda((\CO_{X/G}-q_*(\CM)_+)\otimes(\CO_{X/G}-J)^{\otimes(k+1)})
\end{eqnarray*}
Equality (a) is justified by equality \refeq{veeq}. 
Equality (b) is justified by the adjunction formula. Equality (c) follows 
from the existence of the exact sequence \refeq{exseq}. Equality 
(d) follows from the equality $t\cdot P_k(t)=2^{k+1}-(2-t)^{k+1}$. Equality 
(e) follows again from \refeq{veeq}. Equality (f) is a simple cancellation and so 
is equality (g). Equality (h) follows from the projection formula and the fact that we have 
$q^*(\CJ)\simeq\CL\TW$. Equality (i) follows from the definition of $\lambda(\cdot)$. 
Equality (j) is a simple cancellation and so is equality (k).

Now if we let $k=d$, we obtain by Theorem \ref{DucrotT} 
canonical trivialisations
$$
\lambda((\CO_{X/G}-q_*(\CM)_-)\otimes(\CO_{X/G}-J)^{\otimes(k+1)})\simeq\CO_S
$$
and 
$$
\lambda((\CO_{X/G}-q_*(\CM)_+)\otimes(\CO_{X/G}-J)^{\otimes(k+1)})\simeq\CO_S
$$
and thus a canonical isomorphism
\begin{equation}
\lambda(\iota^*(\CM)\otimes P_d(\CO_{X_G}+N\TW))\simeq\lambda(\CM)^{\otimes 2^{d+1}}.
\label{sieq}
\end{equation}

Note that all the isomorphisms (a),$\dots$, (k) are compatible with any base change to a noetherian scheme.
This follows from that fact that $X\to S$ and $X_G\to S$ are flat, from Lemma \ref{univbc} and from Theorem \ref{DucrotT}.

Now if $X_G$ is not a Cartier divisor let $\wt{X}$ be the blow-up of 
$X$ along $X_G$ and let \mbox{$b:\wt{X}\to X$} be the canonical morphism. 
The scheme $\wt{X}$ is canonically $G$-equivariant since 
the sheaf of ideals of $X_G$ is equivariant. 
 The exceptional divisor $E$ of $\wt{X}$ is isomorphic to the projectivised 
 bundle $\mP(N)$. Since $G$ acts by multiplication by $-1$ on $N$, we see that the 
 action of $G$ is trivial on $E$. Hence $E=\wt{X}_G$ and $\wt{X}_G$ is a Cartier divisor. 
 
 Let $\mu:\wt{X}_G\hookrightarrow \wt{X}$ and 
 $p:\wt{X}_G\to X_G$  be the canonical morphisms. 
 From equality \refeq{sieq}, we obtain
 \begin{eqnarray*}
 &&
 \lambda(\mu^*(b^*(\CM))\otimes P_d(\CO_{\wt{X}_G}+N_{\wt{X}_G/\wt{X}}\TW))
\\&\stackrel{(l)}{\simeq}&
 \lambda(\iota^*(\CM)
 \\&\otimes& 
 \R^\bullet p_*\Big(\CO_{\wt{X}_G}^{\oplus 2^d}+2^{d-1}(\CO_{\wt{X}_G}^{\oplus 2}-(\CO_{\wt{X}_G}+N_{\wt{X}_G/\wt{X}}\TW))+2^{d-2}(\CO_{\wt{X}_G}^{\oplus 2}-
 (\CO_{\wt{X}_G}+N_{\wt{X}_G/\wt{X}}\TW))^{\otimes 2}+\dots
 \\&+&(\CO_{\wt{X}_G}^{\oplus 2}-(\CO_{\wt{X}_G}+N_{\wt{X}_G/\wt{X}}\TW))^{\otimes d}\Big))
 \\&\stackrel{(m)}{\simeq}&
  \lambda(\iota^*(\CM)
 \\&\otimes&
\R^\bullet p_*\Big(\CO_{\wt{X}_G}^{\oplus 2^d}+2^{d-1}(\CO_{\wt{X}_G}-N_{\wt{X}_G/\wt{X}}\TW)+2^{d-2}(\CO_{\wt{X}_G}-N_{\wt{X}_G/\wt{X}}\TW)^{\otimes 2}+\dots
 \\&+&(\CO_{\wt{X}_G}-N_{\wt{X}_G/\wt{X}}\TW)^{\otimes d}\Big))
 \\&\stackrel{(n)}{\simeq}&
 \lambda(\iota^*(\CM)\otimes
 \R^\bullet p_*\Big(
 \sum_{i=0}^d\sum_{j=0}^i 2^{d-i}(-1)^j{i\choose j}(N_{\wt{X}_G/\wt{X}}\TW)^{\otimes j}\Big))
 \\&\stackrel{(o)}{\simeq}&
 \lambda(b^*(\CM))^{\otimes 2^{d+1}}\stackrel{(p)}{\simeq}
 \lambda(\CM)^{\otimes 2^{d+1}}
 \end{eqnarray*}
 For equality (l), use the projection formula. Equality (m) is a simple cancellation. 
 Equality (n) follows from the binomial formula. Equality (o) follows from \refeq{sieq}. 
 Equality (p) follows from the projection formula and the fact that $\R^\bullet b_*(\CO_{\wt{X}})=\CO_X$ (see \cite[VI, §4, proof of Prop. 4.1]{FL-RRA} for lack of a better reference). 
 
Now since $\wt{X}_G=\mP(N)$ we have 
$$
\R^\bullet p_*(N_{\wt{X}_G/\wt{X}}\TW^{\otimes j})\simeq\Sym^j(N\TW)
$$
and we obtain
$$
\lambda(\CM)^{\otimes 2^{d+1}}\simeq 
 \lambda\Big(\iota^*(\CM)\otimes
 \sum_{i=0}^d\sum_{j=0}^i 2^{d-i}(-1)^j{i\choose j}\Sym^j(N\TW)\Big).
$$
Using the fact that there is an equivariant isomorphism $$\Sym^j(N\TW)\simeq(\TW)^{\otimes j}\otimes\Sym^j(N)$$ and using equality \refeq{veeq} we finally get
$$
\lambda(\CM)^{\otimes 2^{d+1}}\simeq 
 \lambda\Big(\iota^*(\CM)\otimes
 \sum_{i=0}^d\sum_{j=0}^i 2^{d-i}{i\choose j}\Sym^j(N)\Big).
$$
Note again that conditional on hypothesis (H) this isomorphism is invariant under any base change to a noetherian scheme by Lemma \ref{lemibl} 
and by the fact that it is invariant under any base change to a noetherian scheme when 
$X_G$ is a Cartier divisor.\endProof

\section{Local refinement of the Adams-Riemann-Roch formula}

\label{secarr}

Let now $g:Y\to S$ be a smooth and strongly projective morphism of noetherian 
schemes. We suppose that $2$ is invertible on $S$. We shall write 
$$X:=Y\times_S Y$$ and we shall write $\pi_1:X\to Y$ and $\pi_2:X\to Y$ for the two projections. 
The group scheme $G=\mZ/2\mZ$ acts on $X$ by swapping the coordinates, with fixed point scheme 
the relative diagonal $\Delta$. The diagonal is regularly immersed 
since $f$ is smooth. 
Let $\CL$ be a line bundle on $Y$. The line bundle $\CM=\pi_1^*(\CL)\otimes\pi_2^*(\CL)$ 
is naturally $G$-equivariant and $\CM|_\Delta\simeq\CL^{\otimes 2}$ 
carries the trivial action. Furthermore $N_{\Delta/X}\simeq\Omega_{Y/S}$ by definition.
Also, note that by the Künneth formula, we have a canonical 
isomorphism
$$
\lambda(\CM)\simeq\lambda(\CL)^{\otimes 2}
$$
where $\lambda(\CM)$ is computed using the above equivariant structure on $\CM.$
Thus applying Theorem \ref{invfunc}, we get an isomorphism
\begin{equation}
\lambda(\CL)^{\otimes 2^{2d+2}}\simeq 
 \lambda\Big(\CL^{\otimes 2}\otimes
 \sum_{i=0}^{2d}\sum_{j=0}^i 2^{2d-i}(-1)^j{i\choose j}\Sym^j(\Omega_{Y/S})\Big).
 \end{equation}
and this completes the proof of Theorem \ref{mainth}.

\end{document}